\documentclass[12pt]{amsart}
\usepackage{amssymb,amsmath,amsfonts,latexsym,setspace}
\usepackage{bm}
\newcommand{\newsection}[1]{\setcounter{equation}{0} \section{#1}}
\newcommand{\bea}{\begin{eqnarray}}
\newcommand{\eea}{\end{eqnarray}}

\newcommand{\clb}{\mathcal{B}}

\newcommand{\cld}{\mathcal{D}}
\newcommand{\cle}{\mathcal{E}}
\newcommand{\clf}{\mathcal{F}}

\newcommand{\clh}{\mathcal{H}}

\newcommand{\clm}{\mathcal{M}}

\newcommand{\clo}{\mathcal{O}}

\newcommand{\clq}{\mathcal{Q}}

\newcommand{\cls}{\mathcal{S}}

\newcommand{\z}{\bm{z}}
\newcommand{\w}{\bm{w}}

\newcommand{\raro}{\rightarrow}

\def \qed {\hfill \vrule height6pt width 6pt depth 0pt}
\def\textmatrix#1&#2\\#3&#4\\{\bigl({#1 \atop #3}\ {#2 \atop #4}\bigr)}
\def\dispmatrix#1&#2\\#3&#4\\{\left({#1 \atop #3}\ {#2 \atop #4}\right)}
\newcommand{\be}{\begin{equation}}
\newcommand{\ee}{\end{equation}}
\newcommand{\ben}{\begin{eqnarray*}}
\newcommand{\een}{\end{eqnarray*}}

\newcommand{\NI}{\noindent}

\newcommand{\bi}{\begin{itemize}}
\newcommand{\ei}{\end{itemize}}

\newtheorem{Theorem}{\sc Theorem}[section]
\newtheorem{Lemma}[Theorem]{\sc Lemma}
\newtheorem{Proposition}[Theorem]{\sc Proposition}
\newtheorem{Corollary}[Theorem]{\sc Corollary}
\newtheorem{Definition}[Theorem]{\sc Definition}
\newtheorem{Example}[Theorem]{\sc Example}
\newtheorem{Remark}[Theorem]{\sc Remark}
\newtheorem{Note}[Theorem]{\sc Note}
\newtheorem{Question}{\sc Question}
\newtheorem{ass}[Theorem]{\sc Assumption}
\newcommand{\bt}{\begin{Theorem}}
\def\beginlem{\begin{Lemma}}
\def\beginprop{\begin{Proposition}}
\def\begincor{\begin{Corollary}}
\def\begindef{\begin{Definition}}
\def\beginexamp{\begin{Example}}
\def\beginrem{\begin{Remark}}
\def\beginq{\begin{Question}}
\def\beginass{\begin{ass}}
\def\beginnote{\begin{Note}}
\newcommand{\et}{\end{Theorem}}
\def\endlem{\end{Lemma}}
\def\endprop{\end{Proposition}}
\def\endcor{\end{Corollary}}
\def\enddef{\end{Definition}}
\def\endexamp{\end{Example}}
\def\endrem{\end{Remark}}
\def\endq{\end{Question}}
\def\endass{\end{ass}}
\def\endnote{\end{Note}}

\begin{document}

\title[An Invariant Subspace Theorem in several variables - II]{An Invariant Subspace Theorem and Invariant Subspaces of Analytic Reproducing Kernel Hilbert Spaces - II}

\author[Jaydeb Sarkar]{Jaydeb Sarkar}
\address{Indian Statistical Institute, Statistics and Mathematics Unit, 8th Mile, Mysore Road, Bangalore, 560059, India}
\email{jay@isibang.ac.in, jaydeb@gmail.com}




\begin{abstract}
This paper is a follow-up contribution to our work \cite{JS1} where
we discussed some invariant subspace results for contractions on
Hilbert spaces. Here we extend the results of \cite{JS1} to the
context of $n$-tuples of bounded linear operators on Hilbert spaces.
Let $T = (T_1, \ldots, T_n)$ be a pure commuting co-spherically
contractive $n$-tuple of operators on a Hilbert space $\clh$ and
$\cls$ be a non-trivial closed subspace of $\clh$.  One of our main
results states that: $\cls$ is a joint $T$-invariant subspace if and
only if there exists a partially isometric operator $\Pi \in
\clb(H^2_n(\cle), \clh)$ such that $\cls = \Pi H^2_n(\cle)$, where
$H^2_n$ is the Drury-Arveson space and $\cle$ is a coefficient
Hilbert space and $T_i \Pi = \Pi M_{z_i}$, $i = 1, \ldots, n$. In
particular, it follows that a shift invariant subspace of a ``nice''
reproducing kernel Hilbert space over the unit ball in
$\mathbb{C}^n$ is the range of a ``multiplier'' with closed range.
Our work addresses the case of joint shift invariant subspaces of
the Hardy space and the weighted Bergman spaces over the unit ball
in $\mathbb{C}^n$.
\end{abstract}

\maketitle

\newsection{Introduction}

Let $T$ be a bounded linear operator on a separable Hilbert space
$\clh$. Furthermore, assume that $T$ is a contraction (that is,
$\|Tf\| \leq \|f\|$ for all $f \in \clh)$ and $T^{*m} \raro 0$ as $m
\raro 0$, in the strong operator topology. Examples of such
$C_{\cdot 0}$-contractions include the multiplication operator $M_z$
on $H^2_{\cle}(\mathbb{D})$, where $\cle$ is a separable Hilbert
space and $H^2_{\cle}(\mathbb{D})$ is the $\cle$-valued Hardy space
over the unit disk $\mathbb{D} = \{z \in \mathbb{C}: |z| <1 \}$.

One of the cornerstones of operator theory and function theory is
that a non-trivial closed $M_z$-invariant subspace of
$H^2_{\cle}(\mathbb{D})$ is the range of a partially isometric
multiplier \cite{RR}. More precisely, let $\cls$ be a non-trivial
closed subspace of $H^2_{\cle}(\mathbb{D})$. Then $\cls$ is a
$M_z$-invariant subspace of $H^2_{\cle}(\mathbb{D})$ if and only if
there exists a Hilbert space $\clf$ and a partially isometric
multiplier $M_{\Theta}$ with symbol $\Theta \in H^\infty_{\clb(\clf,
\cle)}(\mathbb{D})$ such that \[\cls = M_{\Theta}
H^2_{\clf}(\mathbb{D}) = \Theta H^2_{\clf}(\mathbb{D}).\] An
equivalent formulation is that
\[P_{\cls} = M_{\Theta} M_{\Theta}^*,\] where $P_{\cls}$ denotes the
orthogonal projection of $H^2_{\cle}(\mathbb{D})$ onto $\cls$. This
is the celebrated Beurling-Lax-Halmos theorem, due to A. Beurling
\cite{B}, P. Lax \cite{L} and P. Halmos \cite{H}.

In previous work \cite{JS1}, we have shown that the invariant
subspaces of a $C_{\cdot 0}$ contraction $T$ on $\clh$ are given by
the image of those partially isometric operators $\Pi :
H^2_{\clf}(\mathbb{D}) \raro \clh$ which intertwines $T$ and $M_z$,
that is, $\Pi M_z = T \Pi$. It also follows that the shift invariant
subspaces of a large class of reproducing kernel Hilbert spaces also
can be parameterized by the set of all partially isometric
multipliers. These results provide a unifying framework for numerous
invariant subspace theorems, including in particular the
Beurling-Lax-Halmos theorem (see \cite{NF}, \cite{JS}, \cite{DP})
and shift-invariant subspace theorem for weighted Bergman spaces
over the unit disc $\mathbb{D}$ (see \cite{ARS}, \cite{BV1},
\cite{BV2}, \cite{MR}, \cite{O}6).

With this motivation, in this paper, we consider a generalization of
the classification of invariant subspaces of $C_{\cdot
0}$-contractions to pure co-spherically contractive $n$-tuple of
commuting bounded linear operators.

Our present approach is an attempt to understand the joint invariant
subspaces of tuples of commuting operators in the study of operator
theory and function theory in several complex variables. The proofs
of the results in this paper exploit systematically the well known
properties of dilation theory and multiplier spaces of reproducing
kernel Hilbert spaces \cite{A}, so they become simple and clear.
Although our methods of proof are similar to those used in
\cite{JS1}, the present paper extends significantly the class of
multiplication operators tuples on reproducing kernel Hilbert spaces
for which classification of invariant subspaces is known to hold
(cf. \cite{MT}, \cite{GRS}). In particular, we obtain a complete
characterization of shift invariant subspaces of the Hardy space,
the Bergman space and the weighted Bergman spaces over the unit ball
in $\mathbb{C}^n$.

Finally it is worth noting that, among other things, Theorem
\ref{MT} address the following basic question: Let $\clh$ denote the
Drury-Arveson space, or the Hardy space, or the Bergman space, or a
weighted Bergman space over $\mathbb{B}^n$ and $\cls$ be a
non-trivial closed joint $(M_{z_1}, \ldots, M_{z_n})$-invariant
subspace of $\clh$. Does there exists a ``multiplier'' with closed
range such that $\cls$ is the range of the multiplier?

However, in the Drury-Arveson space case, this issue was addressed
by McCullough and Trent \cite{MT} (see also \cite{GRS}).

The paper is organized as follows. Section 2 below contains the
background material on commuting tuples of operators on Hilbert
spaces. Section 3 contains a characterization of invariant subspaces
of pure commuting co-spherically contractive tuples, and Section 4
presents the invariant subspace theorem for reproducing kernel
Hilbert spaces.

\vspace{0.2in}

\NI\textsf{List of symbols:}\begin{enumerate}
\item All Hilbert spaces considered in this paper are separable and
over $\mathbb{C}$. We denote the set of natural numbers including
zero by $\mathbb{N}$.

\item Let $\clh$ be a Hilbert space and $\cls$ be a closed
subspace of $\clh$. The orthogonal projection of $\clh$ onto $\cls$
is denoted by $P_{\cls}$.

\item Let $\clh_1, \clh_2$ and $\clh$ be Hilbert spaces. We denote
by $\clb(\clh_1, \clh_2)$ the set of all bounded linear operators
from $\clh_1$ to $\clh_2$ and $\clb(\clh) = \clb(\clh, \clh)$.

\item Let $n \geq 1$ and $n \in \mathbb{N}$. The set of multi-indices
will be denoted by $\mathbb{N}^n$. That is, $\mathbb{N}^n = \{\bm{k}
= (k_1, \ldots, k_n): k_i \in \mathbb{N}\}$.

\item For $\{z_i\}_{i=1}^n \subseteq \mathbb{C}$, we denote $(z_1,
\ldots, z_n) \in \mathbb{C}^n$ by $\z$.

\item For each $\bm{k} \in \mathbb{N}^n$, define $z^{\bm{k}} =
z_1^{k_1} \cdots z_n^{k_n}$.

\item $\mathbb{B}^n = \{ \z \in \mathbb{C}^n : \|\z\|_{\mathbb{C}^n}
< 1\}$.
\end{enumerate}

\section{Preliminaries}

A commuting $n$-tuple ($n \geq 1$) of bounded linear operators $T =
(T_1, \ldots, T_n)$ is said to be {\it co-spherically contractive},
or define a {\it row contraction}, if
$$\|\sum_{i=1}^{n} T_i h_i\|^2 \leq \sum_{i=1}^{n} \|h_i\|^2, \quad
(h_1, \ldots, h_n \in \clh),$$ or, equivalently, if \[\sum_{i=1}^{n}
T_i T_i^* \leq I_{\clh}.\] Define the \textit{defect operator} and
the \textit{defect space} of $T = (T_1, \ldots, T_n)$ as $D =
(I_{\clh} - \sum_{i=1}^n T_i T_i^*)^{\frac{1}{2}} \in \clb(\clh)$
and $\cld = \overline{\mbox{ran}} D$, respectively.

Natural examples of commuting co-spherically contractive tuples are
the multiplication operator tuples $(M_{z_1}, \ldots, M_{z_n})$ on
the Drury-Arveson space \cite{Ar}, the Hardy space, the Bergman
space and the weighted Bergman spaces (see \cite{AE}, \cite{Guo},
\cite{ZZ}, or Proposition \ref{Prop}) all defined over
$\mathbb{B}^n$. Recall that the Drury-Arveson space, denoted by
$H^2_n$, is determined by the kernel function
\[K_1(\z, \w) = (1 - \sum_{i=1}^n z_i \bar{w}_i)^{-1}. \quad \quad
(\z, \w \in \mathbb{B}^n)\]More generally, for each $\lambda \geq
1$, define the  positive definite function $K_\lambda : \mathbb{B}^n
\times \mathbb{B}^n \raro \mathbb{C}$ by
\[K_\lambda(\z, \w) = (1 - \sum_{i=1}^n z_i \bar{w}_i)^{-\lambda}. \quad \quad
(\z, \w \in \mathbb{B}^n)\] Then the Hardy space
$H^2(\mathbb{B}^n)$, the Bergman space $L^2_a(\mathbb{B}^n)$, and
the weighted Bergman spaces $L^2_{a, \alpha}(\mathbb{B}^n)$, with
$\alpha > 0$, are reproducing kernel Hilbert spaces with kernel
$K_{\lambda}$ for $\lambda = n$, $n+1$ and $n + 1 + \alpha$,
respectively.

Let $\cle$ be a Hilbert space. We identify the Hilbert tensor
product $H^2_n \otimes \cle$ with the $\cle$-valued $H^2_n$ space
$H^2_n(\cle)$, or the $\clb(\cle)$-valued reproducing kernel Hilbert
space with kernel function \[(\z, \w) \mapsto K_1(\z, \w) I_{\cle}.
\quad \quad (\z, \w \in \mathbb{B}^n)\] Then
$$H^2_n(\cle) = \{ f \in \clo (\mathbb{B}^n, \cle): f(z) =
\sum_{\bm{k} \in \mathbb{N}^n} a_{\bm{k}} z^{\bm{k}}, a_{\bm{k}} \in
\cle, \|f\|^2 : = \sum_{\bm{k} \in \mathbb{N}^n} \frac{ \|
a_{\bm{k}} \|^2}{\gamma_{\bm{k}}} < \infty \},$$ where
$\gamma_{\bm{k}} = \frac{(k_1 + \cdots + k_n)!}{k_1 ! \cdots k_n!}$
are the multinomial coefficients and $\bm{k} \in \mathbb{N}^n$ (see
\cite{Ar}, \cite{JS}).

Given a co-spherically contractive tuple $T = (T_1, \ldots, T_n)$ on
$\clh$, define the completely positive map $P_T : \clb(\clh)
\rightarrow \clb(\clh)$ by \[P_{T} (X) = \sum_{i=1}^{n} T_i X T^*_i.
\quad \quad \quad (X \in \clb(\clh))\] Note that \[\sum_{i=1}^n T_i
T_i^* = P_{T} (I_{\clh}) \leq I_{\clh},\] implies \[I_{\clh} \geq
P_{T} (I_{\clh}) \geq P^2_{T} (I_{\clh}) \geq \cdots \geq P^m_{T}
(I_{\clh}) \geq \cdots \geq 0.\]It then follows that,
$$P_{\infty}(T) := \mbox{SOT} - \mathop{\lim}_{m \raro \infty}
P_{T}^m (I_{\clh})$$ exists and $0 \leq P_{\infty}(T) \leq
I_{\clh}$. A co-spherically contractive $T$ is said to be
\textit{pure} (cf. \cite{Ar}, \cite{MV}) if
$$P_{\infty}(T) = 0.$$

\section{Invariant subspaces of co-spherically contractive tuples}

It is well known that a pure co-spherically contractive tuple $T$ on
$\clh$ is jointly unitarily equivalent to the compressed
multiplication operator tuple
\[P_{\cls}M_{z}|_{\clq} := (P_{\cls}M_{z_1}|_{\clq}, \ldots, P_{\cls}M_{z_n}|_{\clq}),\]
for some joint $(M_{z_1}^*, \ldots, M_{z_n}^*)$-invariant subspace
$\clq$ of $H^2_n(\cle)$ and a coefficient Hilbert space $\cle$ (cf.
\cite{AE}, \cite{Ar}, \cite{MV}, \cite{JS}). We include a proof of
this fact for the sake of completeness in a relevant form regarding
our purposes.

\begin{Theorem}\label{DA-pure}
Let $T$ be a pure commuting co-spherically contractive tuple on
$\clh$. Then the map $\Pi \in \clb(H^2_n(\cld), \clh)$ defined by
\[\Pi (K_1(\cdot, \bm{w}) \eta) = (I_{\clh} - \sum_{i=1}^n \bar{w}_i
T_i)^{-1} D \eta, \quad \quad(\bm{w} \in \mathbb{B}^n,\,\eta \in
\cld)\]is co-isometric and\[\Pi M_{z_i} = T_i \Pi. \quad \quad (i =
1, \ldots, n)\]Moreover
\[(\Pi^* h)(\w) = D (I_{\clh} - \mathop{\sum}_{i=1}^n w_i
T_i^*)^{-1}h, \quad \quad (\w \in \mathbb{B}^n, \, h \in
\clh),\]and\[H^2_n(\cld) = \overline{\mbox{span}}\{ z^{\bm{k}}
(\Pi^* \clh) : \bm{k} \in \mathbb{N}^n\}.\]
\end{Theorem}

\NI\textsf{Proof.} First, for $\bm{w} \in \mathbb{B}^n$ define $(w_1
I_{\clh}, \ldots, w_n I_{\clh}) \in \clb(\clh^n, \clh)$ by \[(w_1
I_{\clh}, \ldots, w_n I_{\clh}) (h_1, \ldots, h_n) = \sum_{i=1}^n
w_i h_i. \quad \quad (h_1, \ldots, h_n \in \clh)\]Since \[\|(w_1
I_{\clh}, \ldots, w_n I_{\clh})\| = (\sum_{i=1}^n
|w_i|^2)^{\frac{1}{2}} = \|\w\|_{\mathbb{C}^n},\]it follows that
\[\begin{split}\|\mathop{\sum}_{i=1}^n w_i T_i^*\| & = \|(w_1
I_{\clh}, \ldots, w_n I_{\clh})^* (T_1, \ldots, T_n)\|\leq \|(w_1
I_{\clh}, \ldots, w_n I_{\clh})^*\| \|(T_1, \ldots, T_n)\| \\ & =
(\sum_{i=1}^n |w_i|^2)^{\frac{1}{2}} \|\sum_{i=1}^n T_i
T_i^*\|^{\frac{1}{2}} = \|\bm{w}\|_{\mathbb{C}^n} \|\sum_{i=1}^n T_i
T_i^*\|^{\frac{1}{2}} < 1.\end{split}\]We now define $\Pi^* \in
\clb(\clh, H^2_n(\cld))$ by
\[(\Pi^* h)(\bm{z}):= D (I_{\clh} - \mathop{\sum}_{i=1}^n z_i
T_i^*)^{-1} h = \sum_{\bm{k} \in \mathbb{N}^n} (\gamma_{\bm{k}} D
T^{*\bm{k}} h) z^{\bm{k}},\]for $h \in \clh$ and $\bm{z} \in
\mathbb{B}^n$. Since for all $m \geq 1$, \[P_T^m(D^2) =
P_T^m(I_{\clh} - P_T(I_{\clh})) = P_T^m(I_{\clh}) -
P_T^{m+1}(I_{\clh}),\]and since $\{P^m_T(I_{\clh})\}$ forms a
telescoping series, it follows that
\[\begin{split} \|\Pi^* h\|^2 & = \|\sum_{\bm{k} \in \mathbb{N}^n}
(\gamma_{\bm{k}} D T^{*\bm{k}} h) z^{\bm{k}}\|^2 = \sum_{\bm{k} \in
\mathbb{N}^n} \gamma_{\bm{k}}^2 \|D T^{*\bm{k}} h\|^2
\|z^{\bm{k}}\|^2 = \sum_{\bm{k} \in \mathbb{N}^n}\gamma_{\bm{k}}^2
\|D T^{*\bm{k}} h\|^2 \frac{1}{\gamma_{\bm{k}}}\\ & = \sum_{\bm{k}
\in \mathbb{N}^n}\gamma_{\bm{k}} \|D T^{*\bm{k}} h\|^2 =
\sum_{m=0}^\infty \sum_{|\bm{k}| = m}\gamma_{\bm{k}} \|D T^{*\bm{k}}
h\|^2 = \sum_{m=0}^\infty \sum_{|\bm{k}|=m}\gamma_{\bm{k}} \langle
T^{\bm{k}} D^2 T^{*\bm{k}} h, h \rangle
\\ & = \sum_{m=0}^\infty \langle \sum_{|\bm{k}|=m}\gamma_{\bm{k}}  T^{\bm{k}} D^2 T^{*\bm{k}} h, h
\rangle = \sum_{m=0}^\infty \langle P^m_{T}(D^2) h, h \rangle \\ & =
\sum_{m=0}^\infty (\langle P^m_{T}(I_{\clh}) h, h \rangle - \langle
P^{m+1}_{T}(I_{\clh}) h, h \rangle)\\& = \|h\|^2 - \langle \lim_{m
\raro 0} P^m_T(I_{\clh}) h, h \rangle,\end{split}\]for all $h \in
\clh$. Then, by applying $P_\infty(T) = \lim_{l \raro \infty}
P^l_{T}(I_{\clh}) = 0$, we obtain
\[\|\Pi^* h\| = \|h\|. \quad \quad \quad (h \in \clh)\]In other words,
$\Pi$ is a co-isometry. Moreover, for $h \in \clh$ and $\bm{w} \in
\mathbb{B}^n$ and $\eta \in \cld$, we have
\[\begin{split}\langle \Pi (K_1(\cdot, \bm{w}) \eta), h
\rangle_\clh & = \langle K_1(\cdot, \bm{w}) \eta, D (I_{\clh} -
\mathop{\sum}_{i=1}^n w_i T_i^*)^{-1} h\rangle_{H^2_n(\cld)} \\& =
\langle \sum_{\bm{k} \in \mathbb{N}^n} (\gamma_{\bm{k}}
{\bar{w}}^{\bm{k}} \eta) z^{\bm{k}}, \sum_{\bm{k} \in \mathbb{N}^n}
(\gamma_{\bm{k}} D T^{*\bm{k}} h) z^{\bm{k}} \rangle_{H^2_n(\cld)}
\\& = \sum_{\bm{k} \in \mathbb{N}^n} \gamma_{\bm{k}}
\bar{w}^{\bm{k}} \langle T^{\bm{k}} D \eta, h \rangle_{\clh}\\ & =
\langle (I_{\clh} - \sum_{i=1}^n \bar{w}_i T_i)^{-1} D \eta,
h\rangle_{\clh},\end{split}\]which implies that
\[\Pi (K_1(\cdot, \bm{w}) \eta) =  (I_{\clh} - \sum_{i=1}^n
\bar{w}_i T_i)^{-1} D \eta.\]Next, it follows easily that \[\langle
\Pi (z^{\bm{l}} \eta), h\rangle = \langle z^{\bm{l}} \eta,
\sum_{\bm{k} \in \mathbb{N}^n} (\gamma_{\bm{k}} D T^{*\bm{k}} h)
z^{\bm{k}}\rangle = \gamma_{\bm{l}} \|z^{\bm{l}}\|^2 \langle \eta, D
T^{* \bm{l}} h \rangle = \langle T^{\bm{l}} D \eta, h\rangle,
\]where $\eta \in
\cld$ and $\bm{l} \in \mathbb{N}^n$, and hence \[\Pi (z^{\bm{l}}
\eta) = T^{\bm{l}} D \eta. \quad \quad \quad (\bm{l} \in
\mathbb{N}^n, \eta \in \cld)\] Therefore, we have \[\Pi M_z
(z^{\bm{k}} \eta) = \Pi (z^{\bm{k} + e_i} \eta) = T^{\bm{k} + e_i} D
\eta = T_i (T^{\bm{k}} D \eta) = T_i \Pi (z^{\bm{k}} \eta),\]for
$\bm{l} \in \mathbb{N}^n$ and $\eta \in \cld$, proving $\Pi M_{z_i}
= T_i \Pi$ for $i = 1, \ldots, n$.

\NI Finally, since $\overline{\mbox{span}}\{ z^{\bm{k}} (\Pi^* \clh)
: \bm{k} \in \mathbb{N}^n\}$ is a joint $(M_{z_1}, \ldots,
M_{z_n})$-reducing subspace of $H^2_n(\cld)$, we have \[H^2_n(\cle)
= \overline{\mbox{span}}\{ z^{\bm{k}} \Pi^* \clh : \bm{k} \in
\mathbb{N}^n\},\]for some closed subspace $\cle \subseteq \cld$. On
the other hand, \[(I_{H^2_n(\cle)} - \sum_{i=1}^n M_{z_i} M_{z_i}^*)
= P_{\cle},\]yields \[Dh = (\Pi^* h)(0) = P_{\cle} (\Pi^* h) =
(I_{H^2_n(\cle)} - \sum_{i=1}^n M_{z_i} M_{z_i}^*) (\Pi^* h). \quad
\quad (h \in \clh)\]Therefore, $\cld \subseteq \cle$ and hence $\cld
= \cle$. This completes the proof. \qed

Now we present the main theorem of this section.

\begin{Theorem}\label{inv}
Let $T = (T_1, \ldots, T_n)$ be a pure commuting co-spherically
contractive tuple on $\clh$ and $\cls$ be a non-trivial closed
subspace of $\clh$. Then $\cls$ is a joint $T$-invariant subspace of
$\clh$ if and only if there exists a Hilbert space $\cld$ and a
partially isometric operator $\Pi \in \clb(H^2_n(\cle), \clh)$ such
that
\[\Pi M_{z_i} = T_i \Pi,\]and that
\[\cls = \Pi (H^2_n(\cle)).\]
\end{Theorem}

\NI\textsf{Proof.} Let $\cls$ be a non-trivial joint $T$-invariant
closed subspace of $\clh$. We denote by $T|_{\cls} = (T_1|_{\cls},
\ldots, T_n|_{\cls})$ the $n$ tuple of operators on $\cls$. Note
that $T|_{\cls}$ is a commuting tuple and \[\|\sum_{i=1}^n
T_i|_{\cls} h_i\|^2 = \|\sum_{i=1}^n T_i h_i\|^2. \quad \quad (h_1,
\ldots, h_n \in \cls)\] This clearly implies the row contractivity
of $T|_{\cls}$. Using the identity
\[P^m_{T|_{\cls}} (I_{\cls}) = \sum_{|\bm{k}| = m} \gamma_{\bm{k}}
(T|_{\cls})^{\bm{k}} (T|_{\cls})^{* \bm{k}} = \sum_{|\bm{k}| = m}
\gamma_{\bm{k}} P_{\cls} T^{\bm{k}} P_{\cls} T^{* \bm{k}}|_{\cls} =
\sum_{|\bm{k}| = m} \gamma_{\bm{k}} T^{\bm{k}} P_{\cls} T^{*
\bm{k}}|_{\cls},\]for $m \in \mathbb{N}$, we have \[\begin{split}
\langle P^m_{T|_{\cls}} (I_{\cls}) h, h \rangle & = \sum_{|\bm{k}| =
m} \gamma_{\bm{k}} \langle T^{\bm{k}} P_{\cls} T^{* \bm{k}} h, h
\rangle = \sum_{|\bm{k}| = m} \gamma_{\bm{k}} \|P_{\cls} T^{*
\bm{k}} h\|^2 \leq \sum_{|\bm{k}| = m} \gamma_{\bm{k}} \| T^{*
\bm{k}} h\|^2
\\ & = \langle P^m_{T} (I_{\clh}) h, h \rangle. \quad \quad \quad (h
\in \cls)
\end{split}\]From this and the fact that $P^m_{T} (I_{\clh}) \raro
0$, in the strong operator topology, it readily follows that
\[P^m_{T|_{\cls}} (I_{\cls}) \raro 0,\]in the strong operator
topology. By Theorem \ref{DA-pure} applied to the pure
co-spherically contractive tuple $T|_{\cls}$, there exists a Hilbert
space $\cle$ and a co-isometric map $\Pi_{\cls} : H^2_n(\cle) \raro
\cls$ such that
\[\Pi_{{\cls}} M_{z_i} = T_i|_{\cls} \Pi_{{\cls}}. \quad \quad (i = 1, \ldots,
n)\]Now, consider the inclusion map $i_{\cls} : \cls \raro \clh$.
The properties of the inclusion map imply immediately that
$i_{\cls}$ is an isometry and
\[i_{\cls} T_{j}|_{\cls} = T_j i_{\cls}. \quad \quad \quad (j = 1, \ldots,
n)\]Define $\Pi : H^2_n(\cle) \raro \clh$ by \[\Pi = i_{\cls}
\Pi_{{\cls}}.\] It follows that \[\Pi M_{z_j} = i_{\cls} \,
\Pi_{{\cls}} M_{z_j} = i_{\cls} T_j|_{\cls} \Pi_{{\cls}} =  T_j
i_{\cls} \Pi_{{\cls}} = T_j \Pi,\]for $j = 1, \ldots, n$, and \[\Pi
\Pi^* = (i_{\cls} \, \Pi_{{\cls}}) (\Pi_{{\cls}}^* \, i_{\cls}^*) =
i_{\cls} i_{\cls}^* = P_{\cls}.\]Thus $\Pi$ is partially isometric
and $\mbox{ran} \,\Pi = \cls$. This proves the necessary part.

\NI The sufficient part follows easily from the intertwining
property $T_i \Pi = \Pi M_{z_i}$, for all $i = 1, \ldots, n$, and
the fact that $\cls = \Pi(H^2_n(\cle))$. This completes the proof.
\qed

Also, the joint invariant subspaces of pure co-spherically
contractive tuples can be characterized by the following corollary.

\begin{Corollary}\label{c1}
Let $T = (T_1, \ldots, T_n)$ be a pure commuting co-spherically
contractive tuple on $\clh$ and $\cls$ be a non-trivial closed
subspace of $\clh$. Then $\cls$ is a joint $T$-invariant subspace of
$\clh$ if and only if there exists a Hilbert space $\cle$ and a
bounded linear operator $\Pi \in \clb(H^2_n(\cle), \clh)$ such that
$\Pi M_{z_i} = T_i \Pi$, for $i = 1, \ldots, n$, and
\[P_{\cls} = \Pi \Pi^*.\]
\end{Corollary}

\newsection{Invariant subspaces of analtytic Hilbert spaces}

In this section we classify the joint shift invariant subspaces of a
large class of reproducing kernel Hilbert spaces over $\mathbb{B}^n$
by applying the reasonings from the previous section. We begin by
formulating the notion of analytic Hilbert spaces.

Let $K : \mathbb{B}^n \times \mathbb{B}^n \raro \mathbb{C}$ be a
positive definite kernel such that $K(\z, \w)$ is holomorphic in the
$\z$ variables and anti-holomorphic in $\w$ variables. Then the
reproducing kernel Hilbert space $\clh_K$, corresponding to the
kernel function $K$, is a Hilbert space of holomorphic functions on
$\mathbb{B}^n$ (cf. \cite{A}, \cite{JS}). We say that $\clh_K$ is an
\textit{analytic Hilbert space} over $\mathbb{B}^n$ if the following
conditions are satisfied:

\NI (i) the multiplication operators by the coordinate functions,
denoted by $\{M_{z_1}, \ldots, M_{z_n}\}$ and defined by
\[(M_{z_i} f)(\w) = w_i f(\w), \quad\quad \quad \quad (i = 1, \ldots, n)\]are
bounded, and

\NI (ii) the tuple $n$-tuple $(M_{z_1}, \ldots, M_{z_n})$ on
$\clh_K$ is a pure co-spherically contractive tuple on $\clh_K$,
that is,
\[\sum_{i=1}^n M_{z_i} M_{z_i}^* \leq I_{\clh_K},\]and
\[P_{\infty}(M_z) = 0.\]

Common and important examples of analytic Hilbert spaces include the
Drury-Arveson space $H^2_n$. We also give some typical examples of
analytic Hilbert spaces.

\begin{Proposition}\label{Prop}
Let $\lambda \geq 1$ and ${K_{\lambda}} : \mathbb{B}^n \times
\mathbb{B}^n \raro \mathbb{C}$ be the positive definite kernel
defined as\[K_\lambda(\z, \w) = (1 - \sum_{i=1}^n z_i
\bar{w}_i)^{-\lambda}. \quad \quad (\z, \w \in \mathbb{B}^n)\]Then
$\clh_{K_\lambda}$ is analytic.
\end{Proposition}
\NI\textsf{Proof.} If $\lambda = 1$, then $\clh_K = H^2_n$, and
hence the result holds trivially. So assume $\lambda > 1$. Notice
that
\[K_\lambda(\z, \w) = K_1(\z, \w) K_{\lambda -1}(\z, \w), \quad
\quad (\z, \w \in \mathbb{B}^n)\]and $K_{\lambda -1} : \mathbb{B}^n
\times \mathbb{B}^n \raro \mathbb{C}$ is positive a definite kernel
on $\mathbb{B}^n$. By Theorem 2 of \cite{DMS}, there exists a
coefficient Hilbert space $\cle$ and a joint $(M^*_{z_1} \otimes
I_{\cle}, \ldots, M^*_{z_n} \otimes I_{\cle})$-invariant subspace
$\clq$ of $H^2_n(\cle)$ such that
\[(M_{z_1}, \ldots, M_{z_n}) \cong P_{\clq} (M_{z} \otimes
I_{\cle})|_{\clq} := P_{\clq} (M_{z_1} \otimes I_{\cle}, \ldots,
M_{z_n} \otimes I_{\cle})|_{\clq}.\]For $m \in \mathbb{N}$ we have
\[\begin{split}P^m_{P_{\clq} (M_{z} \otimes I_{\cle})|_{\clq}} (I_{\clq})& =
\sum_{|\bm{k}| = m} \gamma_{\bm{k}} \big(P_{\clq} (M_{z} \otimes
I_{\cle})|_{\clq}\big)^{\bm{k}} \big(P_{\clq} (M_{z} \otimes
I_{\cle})|_{\clq}\big)^{*\bm{k}} \\ & = \sum_{|\bm{k}| = m}
\gamma_{\bm{k}} P_{\clq} (M_{z} \otimes I_{\cle})^{\bm{k}} (M_{z}
\otimes I_{\cle})^{*\bm{k}}|_{\clq}.\end{split}\]This and the fact
that $H^2_n$ is analytic readily implies that $P_{\clq} (M_z \otimes
I_{\cle})|_{\clq}$ is co-spherically contractive and
\[P_{\infty}(P_{\clq} (M_{z} \otimes I_{\cle})|_{\clq}) = SOT-\lim_{m
\raro 0} P^m_{P_{\clq} (M_{z} \otimes I_{\cle})|_{\clq}}(I_{\clq}) =
0.\]Therefore $\clh_K$ is analytic. This completes the proof. \qed

Let $\clh_{K_1}$ and $\clh_{K_2}$ be two analytic Hilbert spaces
corresponding to the kernel functions $K_1$ and $K_2$ on
$\mathbb{B}^n$ and $\cle_1$ and $\cle_2$ be two coefficient Hilbert
spaces. An operator-valued map $\Theta : \mathbb{B}^n \raro
\clb(\cle_1, \cle_2)$ is said to be a \textit{multiplier} from
$\clh_{K_1} \otimes \cle_1$ to $\clh_{K_2} \otimes \cle_2$ if
\[\Theta f \in \clh_{K_2} \otimes \cle_2. \quad \quad (f \in
\clh_{K_1} \otimes \cle_1)\]The set of all multipliers from
$\clh_{K_1} \otimes \cle_1$ to $\clh_{K_2} \otimes \cle_2$ is
denoted by $\clm(\clh_{K_1} \otimes \cle_1, \clh_{K_2} \otimes
\cle_2)$. If $\Theta \in \clm(\clh_{K_1} \otimes \cle_1, \clh_{K_2}
\otimes \cle_2)$, then the multiplication operator $M_{\Theta} :
\clh_{K_1} \otimes \cle_1 \raro \clh_{K_2} \otimes \cle_2$ defined
by \[(M_{\Theta} f)(\w) = (\Theta f)(\w) = \Theta(\w) f(\w), \quad
\quad \quad (f \in \clh_{K_1} \otimes \cle_1, \w \in
\mathbb{B}^n)\]is bounded. This fact follows readily from the closed
graph theorem.

The following provides a characterization of intertwining maps
between analytic Hilbert spaces.

\begin{Proposition}\label{lemma}
Let $\clh_{K_1}$ and $\clh_{K_2}$ be two analytic Hilbert spaces
over $\mathbb{B}^n$ such that \[\bigcap_{i=1}^n \mbox{ker} (M_{z_i}
- w_i I_{\clh_{K_1}})^* = \mathbb{C} K_1(\cdot, \w), \quad \quad (\w
\in \mathbb{B}^n)\]and, let $X \in \clb(\clh_{K_1} \otimes \cle_1,
\clh_{K_2} \otimes \cle_2)$. Then
\[X(M_{z_i} \otimes I_{\cle_1}) = (M_{z_i} \otimes I_{\cle_2}) X,\quad \quad (i = 1, \ldots, n)\]if and
only if $X = M_{\Theta}$ for some $\Theta \in \clm(\clh_{K_1}
\otimes \cle_1, \clh_{K_2} \otimes \cle_2)$.
\end{Proposition}
\NI\textsf{Proof.} Let $X \in \clb(\clh_{K_1} \otimes \cle_1,
\clh_{K_2} \otimes \cle_2)$ and $X(M_{z_i} \otimes I_{\cle_1}) =
(M_{z_i} \otimes I_{\cle_2}) X$, for all $i = 1, \ldots, n$. If $i =
1, \ldots, n$, $\zeta \in \cle_2$ and $\w \in \mathbb{B}^n$, then
\[\begin{split}(M_{z_i} \otimes I_{\cle_1})^* [X^* (K_2(\cdot, \w) \otimes
\zeta)] & = X^* (M_{z_i} \otimes I_{\cle_2})^* (K_2(\cdot, \w)
\otimes \zeta)\\ & = \bar{w}_i [X^* (K_2(\cdot, \w) \otimes
\zeta)].\end{split}\]Thus\[ X^* (K_2(\cdot, \w) \otimes \zeta) \in
\bigcap_{i=1}^n \mbox{ker} \big((M_{z_i} \otimes I_{\cle_1}) -
w_i\big)^*.\] Using this with the fact that \[ \bigcap_{i=1}^n
\mbox{ker} \Big(M_{z_i} - w_i I_{\clh_{K_1}}\Big)^* = \mathbb{C}
K_1(\cdot, \w),\]we have
\[X^* (K_2(\cdot, \w) \otimes \zeta) = K_1(\cdot, \w) \otimes
X(\w) \zeta, \quad \quad \quad (\zeta \in \cle_2)\]for some linear
map $X(\w) : \cle_2 \raro \cle_1$, and for all $\w \in
\mathbb{B}^n$. Moreover,
\[\|X(\w) \zeta\|_{\cle_1} = \frac{1}{\|K_1(\cdot, \w)\|_{\clh_{K_1}}} \|X^* (K_2(\cdot, \w) \otimes
\zeta)\|_{\clh_{K_1} \otimes \cle_1} \leq \frac{\|K_2(\cdot,
\w)\|_{\clh_{K_2}}}{\|K_1(\cdot, \w)\|_{\clh_{K_1}}} \|X\|
\|\zeta\|_{\cle_2},\]for all $\w \in \mathbb{B}^n$ and $\zeta \in
\cle_2$. Therefore $X(\w)$ is bounded and $\Theta(\w) := X(\w)^* \in
\clb(\cle_1, \cle_2)$ for each $\w \in \mathbb{B}^n$. Thus\[X^*
(K_2(\cdot, \w) \otimes \zeta) = K_1(\cdot,\w) \otimes \Theta(\w)^*
\zeta. \quad \quad \quad (\w \in \mathbb{B}^n, \,\zeta \in \cle_2)\]
In order to prove that $\Theta(\w)$ is holomorphic we compute
\[\begin{split}\langle \Theta(\w) \eta, \zeta \rangle_{\cle_2} & = \langle \eta, \Theta(\w)^*
\zeta \rangle_{\cle_1} = \langle K_1(\cdot, 0) \otimes \eta,
K_1(\cdot, \w) \otimes \Theta(\w)^* \zeta \rangle_{\clh_{K_1}
\otimes \cle_1}  \\& = \langle X(K_1(\cdot, 0) \otimes \eta),
K_2(\cdot, \w) \otimes \zeta \rangle_{\clh_{K_2} \otimes \cle_2}.
\quad \quad (\eta \in \cle_1, \, \zeta \in \cle_2)\end{split}\]Since
$\w \mapsto K_2(\cdot, \w)$ is anti-holomorphic, we conclude that
$\w \mapsto \Theta(\w)$ is holomorphic. Hence $\Theta \in
\clm(\clh_{K_1} \otimes \cle_1, \clh_{K_2} \otimes \cle_2)$.

\NI If $\eta \in \cle_1$, $\zeta \in \cle_2$ and $\z, \w \in
\mathbb{B}^n$, then \[\begin{split} \langle X(K_1(\cdot, \w) \otimes
\eta), K_2(\cdot, \z) \otimes \zeta \rangle_{\clh_{K_2} \otimes
\cle_2} & = \langle (K_1(\cdot, \w) \otimes \eta), X^*(K_2(\cdot,
\z) \otimes \zeta) \rangle_{\clh_{K_1} \otimes \cle_1} \\& = \langle
(K_1(\cdot, \w) \otimes \eta), K_1(\cdot, \z) \otimes \Theta(\z)^*
\zeta \rangle_{\clh_{K_1} \otimes \cle_1}\\ & = K_1(\z, \w) \langle
\eta, \Theta(\z)^* \zeta \rangle_{\cle_1}\\ & = K_1(\z, \w) \langle
\Theta (\z) \eta, \zeta \rangle_{\cle_2}\\ & = \langle (M_{\Theta}
(K_1(\cdot, \w) \otimes \eta))(\z), \zeta \rangle_{\cle_2} \\& =
\langle M_{\Theta} (K_1(\cdot, \w) \otimes \eta), K_2(\cdot, \z)
\otimes \zeta \rangle_{\clh_{K_2} \otimes\cle_2}.
\end{split}\]Thus $X = M_{\Theta}$.

\NI Conversely, let $\Theta \in \clm(\clh_{K_1} \otimes \cle_1,
\clh_{K_2} \otimes \cle_2)$. If $f \in \clh_{K_1} \otimes \cle_1$
and $\w \in \mathbb{B}^n$, then \[(z_i \Theta f)(\w) = w_i
\Theta(\w) f(\w) = \Theta(\w) w_i f(\w) = (\Theta z_i f)(\w),\]for
all $i = 1, \ldots, n$. This completes the proof. \qed

The following corollary is a straightforward consequence of
Proposition \ref{lemma} and the fact that, for
$H^2_n$,\[\bigcap_{i=1}^n \mbox{ker} (M_{z_i} - w_i I_{H^2_n})^* =
\mathbb{C} K_1(\cdot, \w). \quad \quad
 (\w \in \mathbb{B}^n)\]

\begin{Corollary}\label{lemma1}
Let $\clh_K$ be an analytic Hilbert space over $\mathbb{B}^n$ and
$\cle$ and $\cle_*$ be two coefficient Hilbert spaces. Let $X$ be in
$\clb(H^2_n \otimes \cle, \clh_{K} \otimes \cle_*)$. Then
\[X(M_{z_i} \otimes I_{\cle_1}) = (M_{z_i} \otimes I_{\cle_2}) X,\quad \quad (i = 1, \ldots, n)\]if and
only if $X = M_{\Theta}$ for some $\Theta \in \clm(H^2_n \otimes
\cle, \clh_{K} \otimes \cle_*)$.
\end{Corollary}

From the previous corollary and Theorem \ref{inv} we readily obtain
the main result of this section.

\begin{Theorem}\label{MT}
Let $\clh_K$ be an analytic Hilbert space over $\mathbb{B}^n$ and
$\cle_*$ be a coefficient Hilbert space. Let $\cls$ be a non-trivial
closed subspace of $\clh_K \otimes \cle_*$. Then $\cls$ is a joint
$(M_{z_1} \otimes I_{\cle_*}, \ldots, M_{z_n} \otimes
I_{\cle_*})$-invariant subspace of $\clh_K \otimes \cle_*$ if and
only if there exists a Hilbert space $\cle$ and a partially
isometric multiplier $\Theta \in \clm(H^2_n \otimes \cle, \clh_{K}
\otimes \cle_*)$ such that
\[\cls = \Theta H^2_n(\cle),\]or equivalently, \[P_{\cls} = M_{\Theta} M_{\Theta}^*.\]
\end{Theorem}

As a particular case of this theorem, we recover the following
results of McCullough and Trent on a generalization of the
Beurling-Lax Halmos theorem in the context of shift invariant
subspaces of vector-valued Drury-Arveson space \cite{MT} (see also
\cite{GRS}).

\begin{Corollary}\label{BLH-MT}
Let $\cle_*$ be a Hilbert space and $\cls$ be a non-trivial closed
subspace of $H^2_n \otimes \cle_*$. Then $\cls$ is a joint $(M_{z_1}
\otimes I_{\cle_*}, \ldots, M_{z_n} \otimes I_{\cle_*})$-invariant
subspace of $H^2_n \otimes \cle_*$ if and only if there exists a
Hilbert space $\cle$ and a partially isometric multiplier $\Theta
\in \clm(H^2_n \otimes \cle, H^2_n \otimes \cle_*)$ such that $\cls
= \Theta (H^2_n \otimes \cle)$.
\end{Corollary}

The classification result, Theorem \ref{MT}, is completely new even
for the case of Hardy space and for the case of weighted Bergman
spaces over $\mathbb{B}^n$.

\end{document}